\documentclass[12pt,reqno]{amsart}

\usepackage{amsmath,amssymb,amsthm,mathtools}
\usepackage[margin=1in]{geometry}
\usepackage{microtype}
\usepackage{enumitem,needspace}
\usepackage[colorlinks=true,linkcolor=blue,citecolor=blue,urlcolor=blue]{hyperref}
\hypersetup{pdftitle={Average character degree below 14/3}}

\newcommand{\Irr}{\operatorname{Irr}}
\newcommand{\acd}{\operatorname{acd}}
\newcommand{\cd}{\operatorname{cd}}

\newcommand{\PSL}{\operatorname{PSL}}
\newcommand{\PGL}{\operatorname{PGL}}
\newcommand{\GL}{\operatorname{GL}}
\newcommand{\SL}{\operatorname{SL}}
\newcommand{\Sp}{\operatorname{Sp}}
\newcommand{\PSp}{\operatorname{PSp}}
\newcommand{\OO}{\operatorname{O}}

\numberwithin{equation}{section}

\newtheorem{theorem}{Theorem}[section]
\newtheorem{lemma}[theorem]{Lemma}
\newtheorem{proposition}[theorem]{Proposition}

\theoremstyle{remark}

\title{Average character degree below $14/3$}

\author[D. Y. Ma]{Dingying Ma}
\address{College of Mathematics and Statistics, Yili Normal University, Yining 835000, China}
\email{mady0229@hotmail.com}

\author[Y. Yang]{Yong Yang*}
\address{Department of Mathematics, Texas State University, 601 University Drive, San Marcos, TX 78666}
\email{yang@txstate.edu}


\begin{document}

\begin{abstract}
We prove that a finite group with average irreducible character degree less than $14/3$ is solvable if it has no chief factor isomorphic to $A_5$. This answers a question of Aziziheris and Eivazzadeh.
\end{abstract}

\maketitle

\section{Introduction}

For a finite group $G$, let $\Irr(G)$ be the set of irreducible complex characters of $G$, and let $\cd(G)=\{\chi(1):\chi\in\Irr(G)\}$. The average character degree of $G$ is
\[
   \acd(G)=\frac{1}{|\Irr(G)|}\sum_{\chi\in\Irr(G)}\chi(1).
\]
It was proved in \cite[Theorem~A]{MoretoNguyen} that $\acd(G)<16/5$ implies that $G$ is solvable. This bound is best possible, since $\acd(A_5)=16/5$. In \cite[Theorem~2]{AziziherisEivazzadeh}, it was shown that a finite nonsolvable group with a nonabelian minimal normal subgroup not isomorphic to $A_5$ satisfies
\[
   \acd(G)\ge \acd(\PSL(2,7))=\frac{14}{3}.
\]
The same paper proves that a finite group $G$ with $\acd(G)<14/3$ is solvable if $G$ has no chief factor isomorphic to $A_5$ and $4\notin\cd(G)$ \cite[Corollary~4]{AziziherisEivazzadeh}. The authors ask whether the hypothesis $4\notin\cd(G)$ can be removed. We answer this question affirmatively.

\begin{theorem}\label{thm:main}
Let $G$ be a finite group with no chief factor isomorphic to $A_5$. If
\[
   \acd(G)<\frac{14}{3},
\]
then $G$ is solvable.
\end{theorem}

The bound is best possible, as shown by $\PSL(2,7)$. The proof follows the induction in \cite[proof of Corollary~4]{AziziherisEivazzadeh}, with arguments for irreducible characters of degrees $2$, $3$, and $4$. In degree $4$, we use the classification of finite linear groups and estimates for projective character degrees. The hypothesis allows chief factors isomorphic to $A_5^2$, which occur in both the imprimitive and tensor-product cases.

Throughout, all groups are finite and all characters are complex. We write $C_n$ for the cyclic group of order $n$ and $\OO_\infty(G)$ for the largest solvable normal subgroup of $G$. If $N\lhd G$, then $\Irr(G\mid N)$ denotes the set of irreducible characters of $G$ whose kernels do not contain $N$. If $\lambda\in\Irr(N)$, then $\Irr(G\mid\lambda)=\Irr(\lambda^G)$ denotes the set of irreducible characters of $G$ lying over $\lambda$, and $\acd(\lambda^G)$ denotes their average degree.

\section{Linear groups of small degree}

\begin{lemma}\label{lem:degree4classification}
Let $G\le\GL_4(\mathbb C)$ be finite, nonsolvable, and irreducible, and assume that $G$ has no chief factor isomorphic to $A_5$. Let $\pi:\GL_4(\mathbb C)\to\PGL_4(\mathbb C)$ be the natural map. Then one of the following holds.
\begin{enumerate}[label=\textup{(\roman*)}]
\item The representation is imprimitive with two blocks of dimension $2$. If $L$ is the kernel of the action on these blocks, then its blockwise projective image is $A_5\times A_5$, and $G/L\cong C_2$ interchanges the two factors.
\item The representation is primitive of tensor-product type, and
\[
   \pi(G)\cong(A_5\times A_5)\rtimes C_2,
\]
where the involution interchanges the two factors.
\item The representation is primitive and
\[
   \pi(G)\cong\PSL(2,7),\ A_6,\ S_6,\ A_7,\ \text{or }\PSp_4(3).
\]
\item The representation is primitive, and $\pi(G)$ has a normal subgroup $E\cong C_2^4$ such that $\pi(G)/E\cong A_6$ or $S_6$.
\end{enumerate}
\end{lemma}

\begin{proof}
Suppose first that the representation is imprimitive. Four one-dimensional blocks give an abelian diagonal kernel and a quotient in $S_4$, making $G$ solvable. Thus there are two blocks of dimension $2$. Let $L$ be the kernel of the action on the blocks. Irreducibility gives $G/L\cong C_2$, so $L$ is nonsolvable. Its two projective images are conjugate finite subgroups of $\PGL_2(\mathbb C)$. Since the kernel of the joint projective map consists of blockwise scalar matrices, it is abelian. The classification of finite subgroups of $\PGL_2(\mathbb C)$ therefore forces both images to be $A_5$.

The joint image is a subdirect subgroup of $A_5\times A_5$. By Goursat's lemma, it is either a diagonal subgroup or the full direct product. In the diagonal case, the quotient of $L$ by the blockwise scalar kernel is a chief factor of $G$ isomorphic to $A_5$, contrary to the hypothesis. Thus the image is $A_5\times A_5$, and an element outside $L$ interchanges its two factors. This proves (i).

Now suppose that the representation is primitive. The possibilities for $\pi(G)$ are listed in \cite[Appendix~A, diagrams~(A.1)--(A.3)]{CheltsovShramov}. In diagram~(A.1), which lists the groups preserving a nonsingular quadric, the only nonsolvable group without an $A_5$ chief factor is $(A_5\times A_5)\rtimes C_2$ with the factor interchange. In particular, $A_5\times A_5$ has an $A_5$ chief factor, since each direct factor is normal. The action on the quadric gives the tensor decomposition in (ii).

In diagram~(A.2), the groups have a normal subgroup $C_2^4$. After the solvable groups and those with quotient $A_5$ or $S_5$ are excluded, the remaining quotients are $A_6$ and $S_6$, giving (iv). Finally, excluding $A_5$ and $S_5$ from diagram~(A.3) gives exactly the groups in (iii).
\end{proof}

\begin{lemma}\label{lem:structure}
Let $G$ be a finite group with no chief factor isomorphic to $A_5$. Let $M\lhd G$ be minimal subject to being nonsolvable, let $N\le M$ be an abelian minimal normal subgroup of $G$, and assume that $\chi\in\Irr(G\mid N)$ has degree $3$ or $4$. Set $R=\OO_\infty(M)$ and $Q=M/R$. Then $N\le Z(G)$, with $N\cong C_3$ if $\chi(1)=3$ and $N\cong C_2$ if $\chi(1)=4$. Moreover,
\[
   Q\cong A_5^2,\ \PSL(2,7),\ A_6,\ A_7,\ \text{or }\PSp_4(3).
\]
If $M\le T\le G$, then
\[
   [T/R:Q\,C_{T/R}(Q)]\le2.
\]
If $\chi(1)=3$, then $Q\cong\PSL(2,7)$ or $A_6$, and this index equals $1$. When $Q\cong A_5^2$, the only possible nontrivial outer action is the interchange of the two factors.
\end{lemma}

\begin{proof}
Minimality implies that $M=M'$. Let $K=\ker\chi$. Since $N\not\le K$ and $N$ is minimal normal in $G$, we have $N\cap K=1$. Thus $M\cap K$ is a proper normal subgroup of $G$ contained in $M$, so it is solvable and lies in $R$.

Set
\[
   \overline G=G/K,\qquad \overline M=MK/K.
\]
The group $\overline M$ is a nontrivial perfect normal subgroup of $\overline G$, which is therefore nonsolvable. The image of $R$ is the solvable radical of $\overline M$: the inverse image in $M$ of $\OO_\infty(\overline M)$ is solvable because its kernel $M\cap K$ and its quotient are solvable. Consequently,
\begin{equation}\label{eq:radical-quotient}
   M/R\cong\overline M/\OO_\infty(\overline M).
\end{equation}
Also,
\[
   [K,M]\le K\cap M\le R,
\]
so the action of $G$ on $Q$ factors through $\overline G$.

Suppose first that $\chi(1)=3$. An imprimitive irreducible group of degree $3$ has an abelian diagonal kernel and a quotient in $S_3$, so it is solvable. Thus $\overline G$ is primitive. By \cite[Theorem~4.2 and Corollary~4.4]{Hu}, its projective image is isomorphic to $A_5$, $\PSL(2,7)$, or $A_6$. The first possibility gives a quotient of $G$ isomorphic to $A_5$, contrary to the hypothesis. Write $\pi$ for the projective map and $S=\pi(\overline G)$. The group $\pi(\overline M)$ is a nontrivial perfect normal subgroup of the simple group $S$, and hence equals $S$. The kernel is scalar, so \eqref{eq:radical-quotient} gives
\[
   Q\cong S\cong\PSL(2,7)\quad\text{or}\quad A_6.
\]
The abelian normal subgroup $\overline N=NK/K$ is scalar: its common eigenspaces would otherwise form a nontrivial system of imprimitivity. Since $N\le M=M'$, its scalar image has determinant $1$, and therefore lies in $Z(\SL_3(\mathbb C))$. As $N\cap K=1$, this proves that $N\cong C_3$ and $N\le Z(G)$. Finally, the identification $Q\cong S$ respects conjugation by $G$, and $S$ is the projective image of the whole group $\overline G$. Thus every automorphism of $Q$ induced by $G$ is inner. Hence $T/R=Q C_{T/R}(Q)$ whenever $M\le T\le G$.

Assume now that $\chi(1)=4$. The representation affording $\chi$ identifies $\overline G$ with a subgroup of $\GL_4(\mathbb C)$ satisfying Lemma~\ref{lem:degree4classification}. Let $\pi$ denote the natural projective map.

In case~(i), the perfect group $\overline M$ acts trivially on the two blocks. Its blockwise projective image is a nontrivial perfect normal subgroup of $A_5\times A_5$ invariant under the factor interchange. Hence this image is the full direct product. The kernel is central in $\overline M$ and consists of blockwise scalars, so \eqref{eq:radical-quotient} gives $Q\cong A_5^2$. The same argument with the ordinary projective map applies in case~(ii).

In case~(iii), the projective image of $\overline M$ is respectively $\PSL(2,7)$, $A_6$, $A_6$, $A_7$, or $\PSp_4(3)$. Its kernel is scalar, so these are also the possibilities for $Q$. In case~(iv), the image of $\overline M$ in $\pi(\overline G)/E$ is perfect and normal. It is nontrivial because the inverse image of $E$ is solvable, and hence it is $A_6$. Its kernel is solvable, so $Q\cong A_6$.

The subgroup $\overline N=NK/K$ is an abelian minimal normal subgroup of $\overline G$. If $\overline G$ is primitive, the common eigenspaces of $\overline N$ form a system permuted by $\overline G$. Primitivity forces a single eigenspace, so $\overline N$ is scalar.

In case~(i), the image of $\overline N$ in $A_5^2$ is an abelian normal subgroup and hence trivial. Thus $\overline N$ consists of blockwise scalars. The projection of the perfect group $\overline M$ on either block has determinant~$1$. It follows that
\[
   \overline N\le\{(\varepsilon I_2,\delta I_2):\varepsilon,\delta\in\{1,-1\}\}.
\]
Conjugation by an element interchanging the blocks interchanges the two signs. The only nontrivial minimal invariant subgroup is the diagonal subgroup. Hence $\overline N=\langle-I_4\rangle$, so $\overline N$ is scalar in this case as well.

Since $N\cap K=1$, scalarity gives $[N,G]=1$. Also $N\le M=M'$, so the determinant of the representation is trivial on $N$. Its faithful scalar image therefore lies in $Z(\SL_4(\mathbb C))$, of order~$4$. A central minimal normal subgroup has prime order, and hence $N\cong C_2$.

In cases~(i) and~(ii), the subgroup fixing both factors acts by inner automorphisms, while the remaining coset interchanges them. In case~(i), this also follows from the fact that the normalizer of $A_5$ in $\PGL_2(\mathbb C)$ is $A_5$ itself. In case~(iii), the projective group acts on its simple normal subgroup by inner automorphisms, except possibly for the involution in $S_6/A_6$. In case~(iv), the action factors through $\pi(\overline G)/E$, so the induced group of outer automorphisms again has order at most $2$. The required index is therefore at most $2$ for every $T$.
\end{proof}

\section{Projective character degrees}

For a finite group $Q$ and a multiplier class $\alpha\in H^2(Q,\mathbb C^\times)$, let $\Irr_\alpha(Q)$ denote the irreducible projective characters associated with a fixed factor set representing $\alpha$. Their degrees depend only on $\alpha$. Equivalently, these are the irreducible characters of a suitable finite central extension of $Q$ lying over a fixed central character. We use the character-triple correspondence in \cite[Chapter~11]{Isaacs}.

The average degree on a finite disjoint union is the weighted average of the averages on its nonempty subsets.

\begin{lemma}\label{lem:perfectcenter}
Let $L$ be a finite perfect group with cyclic center of order $m$. If $\psi\in\Irr(L)$ is faithful, then
\[
   m\mid \psi(1).
\]
\end{lemma}

\begin{proof}
Let $\rho$ afford $\psi$ and let $d=\psi(1)$. Since $L=L'$, the determinant $\det\rho$ is trivial, so $\rho(L)\le \SL_d(\mathbb C)$. By Schur's lemma, $Z(L)$ acts by scalar matrices. Faithfulness embeds $Z(L)$ into $Z(\SL_d(\mathbb C))$, which has order $d$. Thus $m\mid d$.
\end{proof}

\begin{lemma}\label{lem:index2list}
Let $\mathcal D$ be a nonempty finite multiset of positive integers, and let $\sigma$ be a degree-preserving permutation of $\mathcal D$ with $\sigma^2=1$. Let $F$ be the multiset of degrees fixed by $\sigma$, and let $P$ contain one degree from each two-element $\sigma$-orbit. Set
\[
 A_0=\frac{2\bigl(\sum_{d\in F}d+\sum_{d\in P}d\bigr)}{2|F|+|P|},
 \qquad
 A_1=\frac{2\bigl(\sum_{d\in F}d+2\sum_{d\in P}d\bigr)}{|F|+2|P|}.
\]
Suppose that $\mathcal E$ is another nonempty multiset of positive integers, with a degree-preserving permutation $\tau$ satisfying $\tau^2=1$. For the diagonal action of $(\sigma,\tau)$ on $\mathcal D\times\mathcal E$, the multiset consisting of two entries $de$ for each fixed pair $(d,e)$ and one entry $2de$ for each two-element orbit has average at least
\[
 \min\{A_0,A_1\}.
\]
\end{lemma}

\begin{proof}
The $\tau$-orbits partition $\mathcal E$. A fixed element of degree $b$ gives entries with average $bA_0$. For a two-element orbit of common degree $b$, each fixed element of $\mathcal D$ gives one entry $2bd$, and each two-element orbit of $\mathcal D$ gives two entries $2bd$. These entries have average $bA_1$. Since $b\ge1$, the weighted average over the $\tau$-orbits is at least $\min\{A_0,A_1\}$.
\end{proof}

In the notation of Lemma~\ref{lem:index2list}, let $n=|\mathcal D|$ and $s=\sum_{d\in\mathcal D}d$, and suppose that $\sigma$ has $k$ two-element orbits of common degrees $a_1,\ldots,a_k$. Then
\begin{equation}\label{eq:index2averages}
   A_0=\frac{2s-2\sum_{i=1}^k a_i}{2n-3k},
   \qquad A_1=\frac{2s}{n}.
\end{equation}
In particular, if every repeated degree is at most $3s/(2n)$, then $A_0\ge s/n$.

\begin{lemma}\label{lem:A5square-projective}
Let $Q=A_5\times A_5$, and let $\alpha$ be a multiplier class of $Q$. The average degree on $\Irr_\alpha(Q)$ is greater than $14/3$. The same bound holds after tensoring with the irreducible projective characters of a centralizing group for a fixed multiplier class, and after adjoining the factor-interchanging outer involution if it preserves the multiplier classes.
\end{lemma}

\begin{proof}
The Schur multiplier of $A_5$ has order $2$, and that of $A_5\times A_5$ is $C_2\times C_2$. The characters in $\Irr_\alpha(Q)$ are tensor products of projective characters of the two factors. On each factor, the trivial multiplier class has degrees
\[
   1,3,3,4,5
\]
and the nontrivial multiplier class has degrees
\[
   2,2,4,6.
\]
These are the ordinary degrees of $A_5$ and the faithful degrees of $2.A_5\cong\SL_2(5)$, respectively; see \cite{Atlas}. If both multiplier classes are trivial, the average is $(16/5)^2=256/25$. If exactly one is nontrivial, the average is
\[
   \frac{14}{4}\frac{16}{5}=\frac{56}{5},
\]
and if both are nontrivial, the average is $(14/4)^2=49/4$. All three numbers exceed $14/3$.

Suppose that $\alpha$ is fixed by the factor interchange. Its restrictions to the two factors are then equal. If both restrictions are trivial, the interchange fixes the five diagonal tensor products and has ten two-element orbits; Lemma~\ref{lem:index2list} gives
\[
   A_0=\frac{79}{5},\qquad A_1=\frac{512}{25}.
\]
If both restrictions are nontrivial, the same lemma gives
\[
   A_0=\frac{128}{7},\qquad A_1=\frac{49}{2}.
\]
Since $Q$ is perfect, it has no nontrivial linear characters, so the action on $\Irr_\alpha(Q)$ is determined by the factor interchange. The assertion now follows from Lemma~\ref{lem:index2list}.
\end{proof}

\begin{lemma}\label{lem:A6blocks}
Let $\alpha$ be a multiplier class of $A_6$. The average degree on $\Irr_\alpha(A_6)$ is at least $14/3$. The same bound holds after tensoring with the irreducible projective characters of a centralizing group for a fixed multiplier class, and after adjoining an outer involution that preserves the multiplier classes.
\end{lemma}

\begin{proof}
The Schur multiplier of $A_6$ is cyclic of order $6$. For the trivial class, the degree multiset is
\[
 1,5,5,8,8,9,10.
\]
For a class of order $2$, the corresponding faithful characters of $2.A_6\cong\SL_2(9)$ have degrees
\[
 4,4,8,8,10,10.
\]
For a class of order $3$, the corresponding faithful characters of $3.A_6$ have degrees
\[
 3,3,6,9,15
\]
\cite[proof of Corollary~4]{AziziherisEivazzadeh}; see also \cite{Atlas}. An outer involution induces a degree-preserving permutation whose square is the identity. In each list, every repeated degree is at most three halves of the average. Thus \eqref{eq:index2averages} and Lemma~\ref{lem:index2list} give, respectively, the lower bounds
\[
 \frac{46}{7},\qquad \frac{22}{3},\qquad \frac{36}{5}.
\]
All exceed $14/3$.

For a class of order $6$, the corresponding characters are faithful irreducible characters of $6.A_6$. Their degrees are divisible by $6$, by Lemma~\ref{lem:perfectcenter}, so their tensor products have degree at least $6$. In an overgroup of index $2$, the irreducible characters above these tensor products have the same or twice the degree, and hence also have degree at least $6$.
\end{proof}

\begin{lemma}\label{lem:L27blocks}
Let $Q=\PSL(2,7)$, and let $\alpha$ be a multiplier class of $Q$. The average degree on $\Irr_\alpha(Q)$ is at least $14/3$. The same bound holds after tensoring with the irreducible projective characters of a centralizing group for a fixed multiplier class, and after adjoining an outer involution that preserves the multiplier classes.
\end{lemma}

\begin{proof}
The Schur multiplier of $Q$ has order $2$. For the trivial class the degrees are
\[
 1,3,3,6,7,8,
\]
with average $14/3$. For the nontrivial class the faithful degrees of $\SL_2(7)$ are
\[
 4,4,6,6,8,
\]
with average $28/5$; see \cite{Atlas}. In either list an outer involution can only interchange equal-degree characters. Lemma~\ref{lem:index2list} gives the lower bounds $14/3$ and $28/5$, respectively.
\end{proof}

\begin{lemma}\label{lem:A7blocks}
Let $Q=A_7$, and let $\alpha$ be a multiplier class of $Q$. The average degree on $\Irr_\alpha(Q)$ is greater than $14/3$. The same bound holds after tensoring with the irreducible projective characters of a centralizing group for a fixed multiplier class, and after adjoining an outer involution that preserves the multiplier classes.
\end{lemma}

\begin{proof}
The Schur multiplier of $A_7$ is cyclic of order $6$. For the trivial class the ordinary degrees are
\[
 1,6,10,10,14,14,15,21,35,
\]
whose average is $14$. For the class of order $2$, the faithful irreducible characters of $2.A_7$ have degrees
\[
 4,4,14,14,20,20,36,
\]
with average $16$; see \cite{Atlas}. Lemma~\ref{lem:index2list} gives the same bounds for the tensor products and for the irreducible characters above them in an overgroup of index $2$.

For a class of order $3$, Lemma~\ref{lem:perfectcenter} shows that every faithful degree is divisible by $3$. Degree $3$ cannot occur: otherwise $A_7$ would embed in $\PGL_3(\mathbb C)$, contrary to the degree-three classification, in which the only nonabelian simple possibilities are $A_5$, $A_6$, and $\PSL(2,7)$; see \cite[Corollary~4.4]{Hu}. Hence every degree is at least $6$. For a class of order $6$, Lemma~\ref{lem:perfectcenter} again gives that every faithful degree is divisible by $6$. The tensor products therefore have degree at least $6$. In an overgroup of index $2$, the irreducible characters above them have the same or twice the degree, and hence also have degree at least $6$.
\end{proof}

\begin{lemma}\label{lem:PSpblocks}
Let $Q=\PSp_4(3)\cong U_4(2)$, and let $\alpha$ be a multiplier class of $Q$. The average degree on $\Irr_\alpha(Q)$ is greater than $14/3$. The same bound holds after tensoring with the irreducible projective characters of a centralizing group for a fixed multiplier class, and after adjoining an outer involution that preserves the multiplier classes.
\end{lemma}

\begin{proof}
The Schur multiplier of $Q$ has order $2$. The ordinary character degrees are
\[
 1,5,5,6,10,10,15,15,20,24,30,30,30,40,40,45,45,60,64,81,
\]
while the faithful irreducible characters of $2.Q\cong\Sp_4(3)$ have degrees
\[
 4,4,20,20,20,20,20,36,36,60,60,60,64,80.
\]
These degree lists are given in \cite{Atlas}. For the ordinary list, $n=20$ and $s=576$. The minimum of $A_0$ in \eqref{eq:index2averages} occurs when the two characters of degree $45$ are interchanged and every other character is fixed, giving
\[
   A_0=\frac{2\cdot576-2\cdot45}{2\cdot20-3}
       =\frac{1062}{37}.
\]
For the faithful list, $n=14$ and $s=504$. The minimum occurs when two characters of degree $60$ are interchanged and every other character is fixed, giving
\[
   A_0=\frac{2\cdot504-2\cdot60}{2\cdot14-3}
       =\frac{888}{25}.
\]
Indeed, for $c=1062/37$ or $c=888/25$, respectively, the inequality $A_0\ge c$ is equivalent to
\[
   2s-2cn+\sum_{i=1}^k(3c-2a_i)\ge0.
\]
The only repeated degree exceeding $3c/2$ is $45$ in the first list and $60$ in the second, and each permits at most one two-element orbit. The last expression therefore has minimum $0$, attained by the interchange described above. In both cases $A_1=2s/n>c$. Lemma~\ref{lem:index2list} gives the stated bounds, both of which exceed $14/3$.
\end{proof}

\begin{lemma}\label{lem:projective-reduction}
Let $M\lhd T$, let $R=\OO_\infty(M)$, and set $Q=M/R$. Assume that $Q$ is perfect and centerless, and let $\eta\in\Irr(R)$ be $T$-invariant. Set
\[
   H=T/R,\qquad C=C_H(Q),\qquad H_0=QC.
\]
Then $H_0=Q\times C$. Let $\alpha$ be a factor set on $H$ associated with the character triple $(T,R,\eta)$. Its restriction to $H_0$ is cohomologous to a product of factor sets on $Q$ and $C$. For this product factor set, the irreducible projective characters of $H_0$ are precisely the tensor products of the corresponding projective characters of $Q$ and $C$.

Suppose in addition that $|H:H_0|\le2$ and that $Q$ and the induced outer action are as in one of Lemmas~\ref{lem:A5square-projective}, \ref{lem:A6blocks}, \ref{lem:L27blocks}, \ref{lem:A7blocks}, or~\ref{lem:PSpblocks}. Then
\[
   \acd(\eta^T)\ge\eta(1)\frac{14}{3}.
\]
\end{lemma}

\begin{proof}
Since $Q$ is centerless, $Q\cap C=1$, so $H_0=Q\times C$. By the character-triple correspondence, the characters of $T$ above $\eta$ have degrees $\eta(1)$ times the degrees of the irreducible $\alpha$-projective characters of $H$; see \cite[Chapter~11]{Isaacs}.

For $q\in Q$ and $c\in C$, the scalar
\[
   b(q,c)=\frac{\alpha((q,1),(1,c))}{\alpha((1,c),(q,1))}
\]
is multiplicative in each variable. For fixed $c$, it is a linear character of the perfect group $Q$ and hence is trivial. We may therefore choose $\alpha$ so that its restriction to $Q\times C$ is the product of factor sets on $Q$ and $C$. The irreducible projective characters of $H_0$ are tensor products, and their average degree is the product of the two averages. The average degree for the factor $C$ is at least $1$, which proves the assertion when $H=H_0$.

Now suppose that $|H:H_0|=2$. We may choose $\alpha$ with values in a finite cyclic subgroup of $\mathbb C^\times$. Let $\widehat H$ be the corresponding central extension, let $Z$ be its central kernel, and let $\zeta\in\Irr(Z)$ be the faithful character corresponding to $\alpha$. Let $\widehat Q$, $\widehat C$, and $\widehat H_0$ be the inverse images of $Q$, $C$, and $H_0$. Since $b(q,c)=1$ for all $q\in Q$ and $c\in C$, we have
\[
   [\widehat Q,\widehat C]=1,\qquad
   \widehat Q\cap\widehat C=Z,\qquad
   \widehat H_0=\widehat Q\widehat C.
\]
Thus the irreducible characters of $\widehat H_0$ over $\zeta$ correspond to pairs from $\Irr(\widehat Q\mid\zeta)$ and $\Irr(\widehat C\mid\zeta)$, and their degrees are products.

Both $\widehat Q$ and $\widehat C$ are normal in $\widehat H$. Conjugation by an element outside $\widehat H_0$ therefore acts separately on these two character sets. Its square lies in $\widehat H_0=\widehat Q\widehat C$ and acts by an inner automorphism on each factor. Consequently the two induced permutations preserve degrees and have square one.

Clifford theory for the cyclic quotient $\widehat H/\widehat H_0$ gives two extensions of each invariant irreducible character and one induced irreducible character from each two-element orbit. These characters lie over $\zeta$. Lemma~\ref{lem:index2list} and the hypotheses on $Q$ show that their average degree is at least $14/3$. The character-triple correspondence gives the required inequality.
\end{proof}

\Needspace{10\baselineskip}
\begin{proposition}\label{prop:quotient}
Let $G$ be a finite group with no chief factor isomorphic to $A_5$. Let $M\lhd G$ be minimal subject to being nonsolvable, and let $N\le M$ be an abelian minimal normal subgroup of $G$. Assume that
\[
   \acd(G/N)\ge\frac{14}{3}
\]
and that there exists $\chi\in\Irr(G\mid N)$ with $\chi(1)\in\{3,4\}$. Then
\[
   \acd(G)\ge\frac{14}{3}.
\]
\end{proposition}

\begin{proof}
Let $R=\OO_\infty(M)$. By Lemma~\ref{lem:structure}, $N\le Z(G)$ and $N$ has order $2$ or $3$. Also $N\le R$, since $N$ is abelian and normal in $M$. Thus
\[
   \Irr(G)=\Irr(G/N)\sqcup
   \bigsqcup_{\substack{\lambda\in\Irr(N)\\\lambda\ne1_N}}
   \Irr(G\mid\lambda).
\]
It is therefore enough to prove that $\acd(\lambda^G)\ge14/3$ for each nontrivial $\lambda\in\Irr(N)$. Fix such a character $\lambda$.

The sets $\Irr(G\mid\eta)$, with $\eta$ ranging over representatives of the $G$-orbits on $\Irr(R\mid\lambda)$, form a partition of $\Irr(G\mid\lambda)$. Let $\eta\in\Irr(R\mid\lambda)$ and let $T=I_G(\eta)=\{g\in G:\eta^g=\eta\}$. If $M\nleq T$, then $R\le M\cap T<M$. A nontrivial perfect group has no proper subgroup of index at most $4$, since the action on cosets would give a nontrivial perfect subgroup of the solvable group $S_4$. Thus
\[
   [G:T]\ge [M:M\cap T]\ge5,
\]
and every irreducible character of $G$ lying over $\eta$ has degree at least $5\eta(1)>14/3$.

Assume now that $M\le T$. By the Clifford correspondence,
\[
   \acd(\eta^G)=[G:T]\acd(\eta^T).
\]
Let $Q=M/R$. By Lemma~\ref{lem:structure},
\[
   [T/R:Q C_{T/R}(Q)]\le2.
\]
Lemma~\ref{lem:projective-reduction}, together with Lemmas~\ref{lem:A5square-projective}--\ref{lem:PSpblocks}, gives
\[
   \acd(\eta^T)\ge\eta(1)\frac{14}{3}\ge\frac{14}{3}.
\]
Hence $\acd(\eta^G)\ge14/3$.

Thus every set in the partition of $\Irr(G\mid\lambda)$ has average degree at least $14/3$, and hence
\[
   \acd(\lambda^G)\ge\frac{14}{3}.
\]
This holds for every nontrivial $\lambda\in\Irr(N)$. Together with the hypothesis $\acd(G/N)\ge14/3$, the displayed decomposition of $\Irr(G)$ gives $\acd(G)\ge14/3$.
\end{proof}

\begin{proof}[Proof of Theorem~\ref{thm:main}]
Assume that the theorem is false and choose a counterexample $G$ of minimal order. Let $M\lhd G$ be minimal subject to being nonsolvable. Then $M=M'$. Choose a minimal normal subgroup $N$ of $G$ contained in $M$.

If $N$ is nonabelian, then \cite[Theorem~2]{AziziherisEivazzadeh} gives $\acd(G)\ge14/3$, because $N\not\cong A_5$. Thus $N$ is abelian. We claim that
\[
   \acd(G/N)<\frac{14}{3}.
\]
Suppose instead that $\acd(G/N)\ge14/3$. The disjoint union
\[
   \Irr(G)=\Irr(G/N)\sqcup\Irr(G\mid N)
\]
and the inequality $\acd(G)<14/3$ imply that the average degree on $\Irr(G\mid N)$ is less than $14/3$. Hence there exists $\chi\in\Irr(G\mid N)$ with $\chi(1)\le4$. Since $N\le M=M'\le G'$, no such character is linear. Therefore
\[
   \chi(1)\in\{2,3,4\}.
\]

Suppose that $\chi(1)=2$ and let $K=\ker\chi$. Since $N\cap K=1$, the image $MK/K$ is a nontrivial perfect normal subgroup of $G/K$. Thus $G/K$ is nonsolvable. The classification of finite subgroups of $\PGL_2(\mathbb C)$ implies that its projective image is $A_5$. This gives a quotient of $G$ isomorphic to $A_5$, contrary to the hypothesis. If $\chi(1)=3$ or $4$, Proposition~\ref{prop:quotient} gives $\acd(G)\ge14/3$, again a contradiction. This proves the claim.

The quotient $G/N$ has no chief factor isomorphic to $A_5$. By the minimality of $G$, it is solvable. Since $N$ is abelian, $G$ is solvable, a contradiction.
\end{proof}

\section*{Acknowledgements}
The author was supported by Simons Foundation Award~\#918096. OpenAI's ChatGPT was used for mathematical exploration and editorial assistance. The author checked the final arguments and is responsible for the contents of the paper.

\section*{Disclosure Statement}
The author reports no conflict of interest.

\section*{Data Availability}
No data were used in this work.

\end{document}